\documentclass{article}
\usepackage{amsmath,amsthm}
\usepackage{amssymb}
\usepackage{hyperref}
\usepackage{authblk}
\usepackage{graphicx}

\newcommand{\E}{\mathrm{E}}
\newcommand{\V}{\mathrm{Var}}
\newcommand{\R}{\mathbb{R}}
\newtheorem{rmk}{Remark}
\newtheorem{prop}{Proposition}
\newtheorem{corol}{Corollary}

\begin{document}

\title{Sensitivity indices for independent groups of variables}

\author[1]{Baptiste Broto}

\author[2]{Fran\c{c}ois Bachoc }

\author[1]{Marine Depecker}

\author[3]{Jean-Marc Martinez}

\affil[1]{CEA, LIST, Universit\'e Paris-Saclay, F-91120, Palaiseau, France}
\affil[2]{Institut de Math\'ematiques de Toulouse, Universit\'e Paul Sabatier, F-31062 Toulouse, France}
\affil[3]{CEA, DEN-STMF, Universit\'e Paris-Saclay, F-91191 Gif-sur-Yvette, France}

\maketitle

\begin{abstract}
In this paper, we study sensitivity indices for independent groups of variables and we look at the particular case of block-additive models. We show in this case that most of the Sobol indices are equal to zero and that Shapley effects can be estimated more efficiently. We then apply this study to Gaussian linear models, and we provide an efficient algorithm to compute the
theoretical sensitivity indices. In numerical experiments, we show that this algorithm compares favourably to other existing methods.
We also use the theoretical results to improve the estimation of the Shapley effects for general models, when the inputs form independent groups of variables.
\end{abstract}

\section{Introduction}
In global sensitivity analysis, we consider a variable of interest $Y$ which is a function of several variables $X_1,...,X_p$. We aim to associate a value to each input variable, that quantifies its impact on the output variable. We call these values "sensitivity indices". The first sensitivity indices for a general framework were introduced by Sobol in 1993 \cite{sobol_sensitivity_1993} and are called "Sobol indices". They are based on the output variance. Althougt they remain very popular, other sensitivity indices have been proposed in the litterature since then (see for example \cite{plischke_global_2013}, \cite{borgonovo_new_2007}, \cite{chastaing_indices_2013}, \cite{fort_new_2013}, \cite{lemaitre_density_2012}...), and we can find in \cite{borgonovo_common_2016} a general framework for defining sensitivity indices based on variances, on densities, or on distributions. These sensitivity indices are very useful in many applications, for example in physics or in the industry. However, many of them suffer from a lack of interpretation when the input variables are dependent.

Recently, Owen used the notion of Shapley value in order to define new variance-based sensitivity indices in \cite{owen_sobol_2014} which are particularly well adapted for dependent input variables, and that are called "Shapley effects". The Shapley value originates from game theory in \cite{shapley_value_1953}. This quantity can be useful in very different fields (see for example \cite{moretti_combining_2008} or \cite{hubert_strategic_2003}). Hence, there is an active strand of litterature focusing on the computation of Shapley values (\cite{colini-baldeschi_variance_2016}, \cite{fatima_linear_2008}...) but they are adapted to a more general framework than sensitivity analysis. Only a few articles focus on this computation in the field of global sensitivity analysis (see \cite{owen_sobol_2014}, \cite{song_shapley_2016}, \cite{owen_shapley_2017}, \cite{iooss_shapley_2017}).

In this paper, we focus on two popular variance-based sensitivity indices (Sobol indices and Shapley effects) when the input variables form independent groups and more particularly when the model is block-additive. We first prove that the sensitivity indices have much simpler expressions in this framework. In particular, we show that many Sobol indices are equal to zero and that the computation of the Shapley effects is much easier in this case than in the general framework. We provide an algorithm for evaluating these sensitivity indices in the linear Gaussian case, which is particularly efficient when the covariance matrix of the input Gaussian distribution is block diagonal. Finally, we suggest a new method to estimate the Shapley effects for general models with independent groups of variables. We show numerically that our theoretical results improve the accuracy of the estimates.

The paper is organized as follows: in Section 2, we recall the definition of Sobol indices and Shapley effects and their properties. In Section 3, we present our main theoretical results about these sensitivity indices with independent groups of input variables. In Section 4, we focus on the linear Gaussian framework and give an algorithm to compute the variance-based sensitivity indices in this case. We then apply our theoretical results about the variance-based sensitivity indices to give an efficient computation when the covariance matrix is block diagonal. In Section 5, we highlight our results by comparing these algorithms with the existing algorithm suggested in \cite{song_shapley_2016} designed to compute the Shapley effects. We give in Section 6 a numerical procedure based on the results of Section 3, for the estimation of the Shapley effects for general models with independent groups of input variables. The proofs are postponed to the appendix.

\section{Background and notations}\label{background}

Let $(\Omega,\mathcal{A},\mathbb{P})$ be a probability space, $E=\times_{i=1}^p E_i$ be endowed with the product $\sigma$-algebra $\mathcal{E}=\otimes_{i=1}^p \mathcal{E}_i$ and $X$ be a random variable from $\Omega$ to $E$ in $L^2$. The $(X_i)_{i\in [1:p]}$ will be the input variables of the model.\\
Let $f:(E,\mathcal{E})\rightarrow (\mathbb{R},\mathcal{B}(\mathbb{R}))$ be in $L^2$ and $Y=f(X)$. We will call $Y$ the output variable. We assume that $\V(Y)\neq 0$.\\
Let $[1:p]$ be the set of all integers between 1 and $p$. If $u\subseteq[1:p]$ and $x\in \R^p$, we will write $x_u:=(x_i)_{i\in u}$.\\
For all $u\subseteq[1:p]$, we write $V_u:=\V(\E(Y|X_u))$. We can now define the Sobol indices (see \cite{saltelli_sensitivity_2000} and \cite{chastaing_indices_2013}) for a group of variables $X_u$ as:
\begin{equation}\label{Sobol}
S_u:=\frac{1}{\V(Y)}\sum_{v\subseteq u}(-1)^{|u|-|v|} V_v,
\end{equation}
where $|u|$ is the cardinal of $u$.
This quantity assesses the impact of the interaction of the inputs $(X_i)_{i\in u}$ on the model output variability.

\begin{rmk}Actually, we can find two different definitions of the Sobol indices in the literature.
As in \cite{iooss_shapley_2017}, we can also define the closed Sobol indices by
\begin{equation}\label{closed_Sobol}
S_u^{cl}:=\frac{V_u}{\V(Y)}.
\end{equation}
The closed Sobol indice $S_u^{cl}$ assesses the total impact of the inputs $(X_i)_{i\in u}$. However, as in \cite{chastaing_indices_2013}, we choose to focus on the Sobol indices defined by Equation \eqref{Sobol}.
\end{rmk}

These sensitivity indices were the first ones introduced in \cite{sobol_sensitivity_1993} for non-linear models. The Sobol indices are well defined and the sum of the Sobol indices over $u\subseteq [1:p]$ is equal to one, even if the inputs are dependent. Moreover, these sensitivity indices are positive when the inputs are independent. However, this is no longer true in the general case, so that the Sobol indices are less amenable to interpretation with dependent input variables.

To overcome this lack of interpretation, many methods have been proposed. The authors of \cite{jacques_sensitivity_2006} suggested to find the independent groups of variables and to consider each group as one input variable. This method quantifies the impact of these groups of variables on the output but can not quantify the impact of each input variable. Besides in \cite{mara_variance-based_2012}, the authors transformed the dependent inputs into independent inputs to calculate the Sobol indices of the new input variables. Another approach to deal with the dependece between the inputs is to consider only the first order Sobol indices $(S_i)_{i\in [1:p]}$ which always remain in $[0,1]$. \cite{xu_uncertainty_2008} decomposed each first order indice $S_i$ into a correlated contribution $S_i^C$ and an uncorrelated contribution $S_i^U$. However this method requires a linear model. Moreover, the drawback of considering the first order Sobol indices is that their sum is no longer equal to one. Mara and a.l. suggested in \cite{mara_non-parametric_2015} an alternative definition for the Sobol indices in the dependent case, allowing them to remain positive.

Each previous variant of Sobol indices has advantages and disadvantages. Nevertheless, Iooss and Prieur exhibited in \cite{iooss_shapley_2017} an analysis supporting the use of Shapley effects as sensitivity indices when the inputs are dependent.
We can define the Shapley effects as in \cite{owen_sobol_2014} for the input variable $X_i$ as:
\begin{equation}\label{Shapley}
\eta_i:=\frac{1}{p\V(Y)}\sum_{u\subseteq -i}  \begin{pmatrix}
p-1\\ |u|
\end{pmatrix} ^{-1}\left(V_{u\cup \{i\}}-V_u \right)
\end{equation}
where $-i$ is the set $[1:p]\setminus \{i\}$.\\
The Shapley effects have interesting properties for global sensitivity analysis. Indeed, there is only one Shapley effect for each variable (contrary to the Sobol indices). Moreover, the sum of all the Shapley effects is equal to $1$ (see \cite{owen_sobol_2014}) and all these values lie in $[0,1]$ even with dependent inputs, which is very convenient for the interpretation of these sensitivity indices.\\

%To conclude this section, we remark that if we replace $Y$ by $Y-\E(Y)$, we do not change the values of $(V_u)_{u\subseteq [1:p]}$, and so of the Sobol indices and the Shapley effects. We then assume without loss of generality that the expectation of the output variable $Y$ is equal to zero.

\section{Variance-based sensitivity indice properties for independent groups of inputs}\label{theoretical}

In this section, we give theoretical results about variance-based sensitivity indices when the inputs form independent groups of variables.

\subsection{Notations for the independent groups}

Let $\mathcal{C}=\{C_1,...,C_k\}$ be a partition of $[1:p]$ such that the groups of random variables $(X_{C_j})_{j\in [1:k]}$ are independent. Let $A_j:=X_{C_j}$. Let us write
$$
Y=f(X_1,...,X_p)=g(A_1,...,A_k).
$$
Is $w\subseteq[1:k]$, we define 
$$
V_w^g:=\V(\E(Y|A_w)).
$$
As the inputs $(A_1,...A_k)$ are independent, the Hoeffding decomposition (see \cite{hoeffding_class_1948} and \cite{vdv00}) of $g$ is given by:
\begin{equation}\label{eq_hoeff}
g(A)=\sum_{w\subseteq [1:k]}g_w(A_w).
\end{equation}
Similarly to \eqref{Sobol}, the Sobol indices of $g$ are given by
\begin{equation}
S_w^g:=\frac{\V(g_w(A_w))}{\V(Y)}=\frac{1}{\V(Y)}\sum_{z\subseteq w}(-1)^{|w|-|z|} V_z^g.
\end{equation}

\begin{rmk}
As the inputs $(A_1,...,A_k)$ are independent, the Sobol indice $S_w^g$ of $g$ is the variance of $g_w$ divided by $\V(Y)$ and so is non-negative. Moreover, we can estimate it without trouble because the quantities $(V_z^g)_{z\subseteq w}$ are simple to estimate (using the Pick-and-Freeze estimators \cite{saltelli_sensitivity_2000}, \cite{gamboa_statistical_2016} for example).
\end{rmk}

We also define 
$$
V_u^{g,w}:=\V(\E(g_w(A_w)|X_u)).
$$
Writing $C_w:=\bigcup_{j\in w} C_j$, we have $V_u^{g,w}=V_{u\cap C_w}^{g,w}$. If $u \subseteq C_w$, let $S_u^{g,w}$ be the Sobol indice of $g_w$:
$$
S_u^{g,w}:=\frac{1}{\V(g_w(A_w))}\sum_{v\subseteq u}(-1)^{|u|-|v|} V_v^{g,w}.
$$
Equivalently, if $i \in C_w$, let $\eta_i^{g,w}$ be the Shapley effect of $g_w$:
$$
\frac{1}{|C_{w}| \V(g_w(A_w)}\sum_{u\subseteq C_w\setminus \{i\}} \begin{pmatrix}
|C_{w}|-1 \\ |u|
\end{pmatrix}^{-1}(V_{u\cup i}^{g,w}-V_{u}^{g,w}).
$$
Finally, for all $i\in [1:p]$, let $j(i)\in [1:k]$ be the index such that $i \in C_{j(i)}$.

\subsection{Main results for general models}

We will study the Sobol indices and the Shapley effects in the case of block-independent variables. First, we show a proposition about the $(V_u)_{u\subseteq [1:p]}$.

\begin{prop}\label{propcentrale}
For all $u \subseteq [1:p]$, we have:
\begin{equation}
V_u=\sum_{w\subseteq[1:k]} V_{u\cap C_w}^{g,w}.
\end{equation}
\end{prop}
 Dividing by the variance of $Y$, we can deduce directly an identical decomposition for the closed Sobol indices defined by Equation \eqref{closed_Sobol}.\\
We then provide a consequence of Proposition \ref{propcentrale} on Sobol indices:

\begin{prop}\label{propSobol}
For all $u\subseteq [1:p]$, we have:
\begin{equation}
S_u=\displaystyle\sum_{\substack{w\subseteq[1:k],\\ u\subseteq C_w}} S_w^g S_u^{g,w}.
\end{equation}
\end{prop}

Proposition \ref{propSobol} improves the interpretation of the Sobol indices. It states that the Sobol indices for the output $Y$ are linear combinations of the Sobol indices when considering the outputs $g_w(A_w)$ and that the weighting coefficients are the $S_w^g$.

This proposition can be beneficial for the estimation of the Sobol indices. We could estimate the coefficients $S_w^g$ by Pick-and-Freeze. If many of them are close to $0$, the corresponding Sobol indices $S_u^{g,w}$ are irrelevant for the total output $Y$, and it is unnecessary to estimate them.

We also provide a consequence of Proposition \ref{propcentrale} on Shapley effects:

\begin{prop}\label{propshapley}
For all $i\in [1:p]$, we have
\begin{equation}
\eta_i=\displaystyle\sum_{\substack{w\subseteq[1:k],\\ \text{s.t. } j(i)\in w}}S_w^g  \eta_i^{g,w}.
\end{equation}
\end{prop}
As for the Sobol indices, Proposition \ref{propshapley} provides the computation of the Shapley effects for the output $Y$ by summing the Shapley effects when considering the outputs $g_w(A_w)$ and multiplying them by the coefficient $S_w^g$.

\subsection{Main results for block-additive models}
In the following, we detail the consequences of Propositions \ref{propSobol} and \ref{propshapley} to the particular case of a block-additive model:
\begin{equation}\label{forme}
Y=\sum_{j=1}^kg_j(A_j),
\end{equation}
i.e. when the functions $(g_w)_{w\subseteq[1:k]}$ of the Hoeffding decomposition are equal to 0 except for $w$ equal to a singleton.

\begin{corol}\label{coroSobol}
If the model is block-additive, for all $u$ such that $u\nsubseteq C_j$ for all $j$, we have $S_u=0$.
\end{corol}
This corollary states that the majority of Sobol indices for block additive models a equal to zero. It remains only $\sum_{j=1}^k 2^{C_j}-1$ unknown non-zero Sobol indices instead of $2^p-1$.

\begin{corol}\label{coroshapley}
For block-additive models, we have
\begin{equation}
\eta_i=S_{j(i)}^g\eta_i^{g,j(i)}.
\end{equation}
\end{corol}
For example, if we apply this corollary in the case where $X_i$ is the only variable in its group, then we have $\eta_i=S_i$.\\
To compute the Shapley effect $\eta_i$ in block additive models, the previous corollary reduces the sum from all the subsets of $[1:p]\setminus \{i\}$ to all the subset of $C_{j(i)} \setminus \{i\}$. Then, the computational gain is the same as in Corollary \ref{coroSobol}.

\section{Explicit computation of sensitivity indices for Gaussian linear models}
In general, it is a really difficult task to estimate the sensitivity indices for dependent inputs. The authors of \cite{song_shapley_2016} suggest an algorithm for the Shapley effects estimation which requires a function generating an i.i.d. sample of $X_{-u}$ conditionally to $X_u=x_u$. This requirement reduces this estimation to restricted theoretical frameworks. In the following, we show an exact computation of the sensitivity indices for Gaussian linear models.
\subsection{Linear Gaussian framework}
In this section, we assume that $X\sim \mathcal{N}(\mu,\Gamma)$, that $\Gamma$ is invertible and that $f:x\longmapsto \beta_0+\beta^T x$, for a fixed $\beta_0\in \R$ and a fixed vector $\beta$. This framework is widely used to model physical phenomena (see for example \cite{kawano_evaluation_2006}, \cite{hammer_approximate_2011}, \cite{rosti_linear_2004}). Indeed, uncertainties are often modelled as Gaussian variables and an unknown function is commonly estimated by its linear approximation. We can assume without loss of generality that $\mu=0$ and $\beta_0=0$. For now, we will not assume that there are different groups of independent variables.
In this framework, the sensitivity indices can be calculated explicitly. First, we write the sensitivity indices with expectations of conditional variances:
\begin{eqnarray}\label{varianceSobol}
S_u:=(-1)^{|u|}\sum_{v\subseteq u} (-1)^{|v|+1} \frac{\E(\V(Y|X_u))}{\V(Y)},\;\forall u \neq \emptyset,
\end{eqnarray}
and
\begin{eqnarray}\label{varianceShapley}
\eta_i:=\frac{1}{p\V(Y)}\sum_{u\subseteq -i}  \begin{pmatrix}
p-1\\ |u|
\end{pmatrix} ^{-1}\left(\E(\V(Y|X_u))-\E(\V(Y|X_{u\cup \{i\}})) \right).
\end{eqnarray}
We will now exploit the Gaussian framework as in \cite{owen_shapley_2017}, using that for all subset $u\subseteq[1:p]$,
\begin{equation}\label{eq}
\V(Y|X_u)=\V(\beta_{-u}^T X_{-u}|X_u)=\beta_{-u}^T(\Gamma_{-u,-u}-\Gamma_{-u,u}\Gamma_{u,u}^{-1}\Gamma_{u,-u})\beta_{-u}
\end{equation}
where $\beta_u:=(\beta_i)_{i\in u}$ and $\Gamma_{u,v}:=(\Gamma_{i,j})_{i\in u, j\in v}$. These conditional variances are constant so they are equal to their expectation. Then, we can use these formulae to compute explicitly these sensitivity indices in the Gaussian linear case.

\begin{rmk}
If the matrix $\Gamma$ is not invertible, there exist subsets $u$ such that $\Gamma_{u,u}$ is not invertible. However, Equation \eqref{eq} still holds if we replace $\Gamma_{u,u}^{-1}$ by the generalized inverse (for symmetric matrices) of $\Gamma_{u,u}$.
\end{rmk}

\begin{rmk}
One can show a similar result when $X$ follows an asymmetric Laplace distribution $AL_p(m,\Gamma)$. However, the conditional variances are not constant in this case and their expectations  must be estimated, for instance by Monte-Carlo.
\end{rmk}

One issue remains though, namely computing numerically the sum in \eqref{varianceShapley}. Indeed, we have to sum over all the subsets of $[1:p]$ which do not contain $i$. We also have to group the subsets $u$ and $u\cup\{i\}$. For this purpose, we suggest to use the following bijective map:
\begin{eqnarray*}
h:
\begin{array}{lll}
\mathcal{P}([1:p])&\longrightarrow &[0:2^{p}-1]\\
u&\longmapsto &\sum_{i\in u} 2^{i-1}.
\end{array}
\end{eqnarray*}
We remark that:
\begin{eqnarray*}
u\subseteq -i 
 \Longleftrightarrow \left\lfloor \frac{h(u)}{2^{i-1}}\right\rfloor \equiv 0\; [mod\;2].
\end{eqnarray*} 
Finally, we can see that if $u\subseteq -i$, then $h(u\cup\{i\})=h(u)+h(\{i\})$. 
Based on this map and Equations \eqref{varianceSobol} and \eqref{varianceShapley}, we suggest an algorithm that we call "LG-Indices" (for Linear Gaussian). This algorithm computes the variance-based sensitivity indices in the linear Gaussian framework.

\underline{\textbf{LG-Indices:}} \textbf{Inputs:} $\beta$, $\Gamma$.
\begin{enumerate}
\item Let $\V(Y)=\beta^T \Gamma \beta$ and let $\V(Y|X)=0$.
\item (Compute the conditional variances) For $j=0,...,2^p-1$, do the following:
\begin{enumerate}
\item Compute $u=h^{-1}(j)$.
\item Compute $V_j:=\V(Y|X_u)$ using \eqref{eq}.
\end{enumerate}
\item Initialise $S=(0,...,0)\in \R^{2^p}$.
\item (Compute the Sobol indices) For $j=0,...,2^p-1$, do the following:
\begin{enumerate}
\item Let $v=h^{-1}(j)$.
\item For all $u$ such that $v\subseteq u$, let:
\begin{equation}
S_{h(u)}=S_{h(u)}+(-1)^{|v|+1}V_{j}
\end{equation}
\end{enumerate}
\item Let $S_0=0$
\item For $j=1,...,2^p-1$, let $u=h^{-1}(j)$ and
\begin{equation}
S_j=\frac{(-1)^{|u|}}{\V(Y)}S_j.
\end{equation}
\item (Compute the Shapley effects) For $i=1,...,p$, do the following:
\begin{enumerate}
\item Initialize $\eta=(0,...,0)\in \R^p$ .
\item For $k=0,...,2^p-1$, do the following:
\begin{enumerate}
\item If $\left\lfloor \frac{k}{2^{i-1}}\right\rfloor \equiv 0\; [mod\;2]$, then update :
\begin{equation}
\eta_i=\eta_i+ \begin{pmatrix}
p-1 \\ |h^{-1}(k)|
\end{pmatrix}^{-1} \left(V_k-V_{k+2^{i-1}})) \right).
\end{equation}
\end{enumerate}
\item Let $\eta_i=\eta_i\slash (p\V(Y))$.
\end{enumerate}
\end{enumerate}
\textbf{Outputs} $(S,\eta)$.
\bigskip

\begin{rmk}
We can use the previous algorithm for any $f$ and any law of $X$ if we can estimate the expectation of conditional variances (or the variance of conditional expectations equivalently).
\end{rmk}

\subsection{Linear Gaussian framework with independent groups of variables}

Despite the analytical formula \eqref{eq}, the computational cost remains an issue when the number of input variables $p$ is too large. Based on an implementation in the R software, LG-Indices provides almost instantaneous results for $p\leq 15$, but becomes impracticable for $p\geq 30$. Indeed, we have to store $2^p$ values, namely the $(\V(Y|X_u))_{u\subseteq[1:p]}$, and this can be a significant issue. Fortunately, when $p$ is large, it can frequently be the case that there are independent groups of random variables. That is, after a permutation of the variables, the covariance matrix $\Gamma$ is a block diagonal matrix. In this case, Corollaries \ref{coroSobol} and \ref{coroshapley} show that this high dimensional computational problem boils down to a collection of lower dimensional problems.

In this framework, we have seen in Corollaries \ref{coroSobol} and \ref{coroshapley} that we only have to calculate the $\sum_{j=1}^k2^{|C_j|}$ values $\{\V(Y|X_u),\; u\subseteq C_j,\;j\in[1:k]\}$ instead of all the $2^p$ values $\{\V(Y|X_u), \; u\subseteq[1:p]\}$. We detail this idea in the algorithm "LG-GroupsIndices".

\underline{\textbf{LG-GroupsIndices:}} \textbf{Inputs:} $\beta$, $\Gamma$.
\begin{enumerate}
\item Let $C_1,...,C_k$ be the independent groups of variables, for example using the function "graph$\_$from$\_$adjacency$\_$matrix" of the R package "igraph" (see \cite{csardi_igraph_2006}).
\item Let $\eta$ be a vector of size $p$.
\item For $j=1,...,k$, do the following:
\begin{enumerate}
\item Let $(\tilde{S},\tilde{\eta})$ be the output of the algorithm LG-Indices with the inputs $\beta_{C_j}$ and $\Gamma_{C_j,C_j}$.
\item Let
\begin{equation*}
S^j=\frac{\beta_{C_j}^T\Gamma_{C,j,C_j} \beta_{C_j}}{\beta_T \Gamma \beta_T}\tilde{S}
\end{equation*}
\item Let 
\begin{equation*}
\eta_{C_j}=\frac{\beta_{C_j}^T\Gamma_{C,j,C_j} \beta_{C_j}}{\beta_T \Gamma \beta_T}\tilde{\eta}.
\end{equation*}
\end{enumerate}
\end{enumerate}
\textbf{Ouputs:} $(S^1,...,S^k,\eta)$.
\bigskip

We have used LG-GroupsIndices for computing Shapley effects on an industrial study in the field of nuclear safety. In this model, the twelve inputs were modelled by a Gaussian vector with two independent groups of six variables. Corollary \ref{coroshapley} enables us to compute the Shapley effects computing only $2^6+2^6=128$ conditional variances instead of $2^{12}=4096$. These results have been presented by Pietro Mosca in the sixteen International Symposium on Reactor Dosimetry.\\
The computational time does not really depend of the coefficients of $\beta$ and $\Gamma$, that is why we prefer to consider simulated toy examples to compare the different algorithms in the next section.

\section{Numerical experiments: Shapley effects computations}
To position our work with respect to the state of art, we compare in this section the algorithms "LG-Indices" and "LG-GroupsIndices" with an existing algorithm designed to compute the Shapley effects for global sensitivity analysis.
\subsection{Random permutations Algorithm: Shapley effects estimation}

To the best of our knowledge, the only existing algorithm designed to compute the Shapley effects for global sensitivity analysis is suggested by Song in \cite{song_shapley_2016}. This algorithm is introduced for a general function $f$ and a general distribution of $X$ and does not focus specifically on the linear Gaussian model. In this general context, the expectation of the conditional variances are estimated by a double Monte-Carlo procedure.

This Shapley effect estimation suggested in \cite{song_shapley_2016} relies on the following formulation of the $\eta_i$:
\begin{equation}\label{formulation}
\eta_i=\frac{1}{p!\V(Y)}\sum_{\sigma \in \mathcal{S}_p}(\E(\V(Y|X_{P_i(\sigma)}))- \E(\V(Y|X_{P_i(\sigma)\cup \{i\})}))),
\end{equation}
where $\mathcal{S}_p$ is the set of all permutations of $[1:p]$ and $P_i(\sigma):=\{\sigma(j),\;j\in [1:i] \}$.
Then, in order to circumvent the $p!$ estimations of expectations of conditional variances, \cite{song_shapley_2016} suggests to generate $m$ ($m<p!$ for large values of $p$) permutations $\sigma_1,...,\sigma_m$ uniformly in $\mathcal{S}_p$ and to let:
\begin{equation}\label{song}
\widehat{\eta}_i=\frac{1}{m\V(Y)}\sum_{j=1}^m(\E(\V(Y|X_{P_i(\sigma_j)}))- \E(\V(Y|X_{P_i(\sigma_j)\cup \{i\})}))).
\end{equation}
It is suggested in \cite{iooss_shapley_2017} to choose $m$ as large as possible and to choose small sample size values for the double Monte-Carlo procedure. This algorithm is already implemented in the R package {\verb sensitivity } and is called "shapleyPermRand". However, as we focus on the linear Gaussian framework, for a fair comparison, we adapt the algorithm suggested in \cite{song_shapley_2016} to this particular framework replacing the estimations of $(\E(\V(Y|X_u)))_{u\subseteq[1:p]}$ by their theoretical values given by \eqref{eq}. We will write "random permutations Algorithm" for this algorithm.

A variant of this algorithm, called "shapleyPermEx", is implemented in the R package "sensitivity". Although, this algorithm is not clearly suggested in \cite{song_shapley_2016}. "shapleyPermEx" differs from "shapleyPermRand" by computing the sum over all the permutations when the former algorithm was only estimating the sum thanks to a Monte-Carlo method. We used the linear Gaussian framework to replace the expectation of conditional variances by their theoretical values, such that the algorithm gives the exact values of the Shapley effects. We will write "exact permutations Algorithm" for this algorithm. This method still remains very costly due to the computation of $(p-1)!$ conditional variances. For example, it spends more than ten minutes computing the Shapley effects for only $p=10$.

\subsection{Shapley effects in the linear Gaussian framework}\label{section_algo1}
From now on, as we focus on the Shapley effects, we do not carry out the steps related to the computation of Sobol indices in the algorithms "LG-indices" and "LG-GroupsIndices".\\

Let us consider a simulated toy example. We generate $\beta$ by a $\mathcal{N}(0,I_p)$ random variable and $\Gamma$ by writing $\Gamma=AA^T$, where the coefficients of $A$ are generated independently with a standard normal distribution.\\

First, we compare "LG-Indices" with exact permutations Algorithm. Both provide the exact Shapley values but with different computational times. Table \ref{permex} provides the computation times in seconds for different values of $p$, the number of input variables.
\begin{table}[ht]
\caption{Computational time (in seconds) for exact permutations Algorithm and LG-Indices for different values $p$.}
\centering
\begin{tabular}{|l|c|c|c|c|c|}\hline
 & $p=6$& $p=7$& $p=8$ & $p=9$ & $p=10$   \\ \hline
exact permutations & 0.11 & 0.78 & 7.31 & 77.6 & 925  \\ \hline
LG-Indices & 0.004 & 0.008& 0.018 & 0.039 & 0.086   \\ \hline
\end{tabular}
\label{permex}
\end{table}
We remark that "LG-indices" is much faster than exact permutations Algorithm. \\

We can also compare "LG-Indices" with random permutations Algorithm. For the latter, we choose $m$, the number of permutations generated in Song, so that the computational time is the same as LG-Indices. Yet, while our algorithm gives the exact Shapley effects, the random permutations Algorithm provides an estimation of them. Hence, the performance of  the latter algorithm is evaluated by computing the coefficients of variation in $\%$, for different values of $p$. We give in Table \ref{Song1} the average of the $p$ coefficients of variations. We recall that the coefficient of variation corresponds to the ratio of the standard deviation over the mean value.
\begin{table}[ht]
\caption{Mean of the coefficients of variation of Shapley effects estimated by random permutations Algorithm for the same computational time as LG-GroupsIndices.}
\centering
\begin{tabular}{|l|c|c|c|c|c|}\hline
 & $p=3$& $p=4$& $p=5$ & $p=6$ & $p=7$   \\ \hline
 $m$ chosen & 10 & 12& 18 &30 & 50  \\ \hline
mean of coefficients of variation &$31\%$ & $26\%$& $22\%$ & $18\%$ & $13\%$   \\ \hline
\end{tabular}
\label{Song1}
\end{table}
We see that the algorithm suggested in \cite{song_shapley_2016} has quite large coefficients of variation when we choose $m$ so that the computational time is the same as our algorithm. However, this variation decreases with the number of inputs $p$. We can explain that by saying that the computational time of LG-Indices is exponential with $p$. So, we can see that the precision of random permutations Algorithm increases with $m$.

\subsection{Shapley effects in the linear Gaussian framework with independent groups of variables}
With independent groups of inputs, to the best of our knowledge, LG-GroupsIndices is the only algorithm which can compute the exact Shapley effects for large values of $p$ (the number of inputs). Indeed, random permutations Algorithm can handle large values of $p$ but always computes estimations of Shapley effects. On the other hand, LG-Indices computes exact Shapley effects but becomes too costly for $p\geq 20$ (the computation time is exponential in $p$).\\

First, we compare the computation time of "LG-Indices" and "LG-GroupsIndices" for low values of $p$ on a toy simulated example as in Section \ref{section_algo1}. We generate $k$ independent groups of $n$ variables. We give these results in Table \ref{algo1_2}.

\begin{table}[ht]
\caption{Computation time (in seconds) for LG-Indices and LG-GroupsIndices for different values of $k$ and $n$.}
\centering
\begin{tabular}{|l|c|c|c|c|}\hline
 & $k=3$& $k=4$& $k=4$  & $k=5$  \\
 & $n=3$&  $n=3$  & $n=4$   & $n=4$     \\ \hline
 LG-Indices & $0.04$ & $0.47$ & $8.45$  & $168.03$  \\ \hline
LG-GroupsIndices &$0.002$ & $0.003$ & $0.004$  & $0.007$ \\ \hline
\end{tabular}
\label{algo1_2}
\end{table}
\bigskip

Now, we compare LG-GroupsIndices with random permutations Algorithm as in \ref{section_algo1}: we choose $m$ so that the computational time is the same and we give the average of the $p$ coefficients of variation of random permutations Algorithm in Table \ref{algoSong_2}.

\begin{table}[ht]
\caption{Mean of the coefficients of variation of Shapley effects estimated by random permutations Algorithm for the same computational time as LG-GroupsIndices.}
\centering
\begin{tabular}{|l|c|c|c|c|c|c|}\hline
 & $k=3$& $k=4$& $k=5$ & $k=6$ & $k=10$ & $k=5$\\
  & $n=3$ & $n=4$ & $n=5$ & $n=6$ & $n=5$ & $n=10$ \\ \hline
 $m$ chosen & 7 & 8& 9 &10 & 5 & 98 \\ \hline
mean of coefficients of variation &$34\%$ & $38\%$& $34\%$ & $36\%$ & $40\%$ & $12\%$  \\ \hline
\end{tabular}
\label{algoSong_2}
\end{table}
Here, the mean of the coefficients of variation remains quite large (around $35\%$) when we chose $k=n$. However, when we choose $k$ larger (resp. lower) than $n$, this variation increases (resp. decreases).

\section{Numerical applications: approximating general models by block-additive functions}

In this section, we give a numerical application of Corollary \ref{coroshapley} to the estimation of the Shapley effects with a general model and with independent groups of input variables.

We assume that:
\begin{itemize}
    \item There is a numerical code $f:\R^p\rightarrow \R$ which is very costly;
    \item It is feasible to sample under the conditional distributions of the input variables (in order to use the algorithm suggested in \cite{song_shapley_2016} to estimate the Shapley effects);
    \item There are independent groups of input variables.
\end{itemize}
Since the numerical code is costly, we need to compute a meta-model $\tilde{f}$ of $f$ to estimate the Shapley effects from a sample $(X_i,f(X_i))_{i\leq N_{samp}}$. We compare two different methods:
\begin{enumerate}
    \item We compute a total meta-model $\tilde{f}_{tot}$ without assumptions on the model. Then, we compute the Shapley effects using the R function "shapleyPermRand" from the R package {\verb sensitivity }.
    \item We compute a block-additive meta-model $\tilde{f}_{add}=\sum_{j=1}^k g_j$ as in Equation \eqref{forme}, where $g_j$ depends only on the group of variable $X_{C_j}$. Then, we estimate the Shapley effects for all the models $g_j$ using the R function "shapleyPermRand". From them, we deduce the Shapley effects of the block-additive model $\tilde{f}_{add}$ using Corollary \ref{coroshapley}.
\end{enumerate}
To compute the two models $\tilde{f}_{tot}$ and $\tilde{f}_{add}$, we suggest to use the functions "kmAdditive" and "predictAdditive" from the R package {\verb fanovaGraph }.

\begin{figure}
    \centering
    \includegraphics[height=8cm,width=10cm]{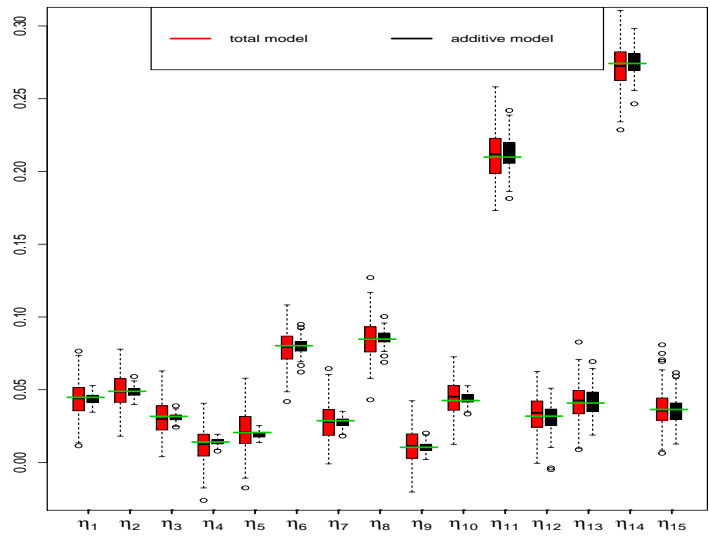}
    \caption{Boxplot of 1000 estimates of the 15 Shapley effects using the two methods suggested: with a total model in red and with an additive model in black. We also represent the true values of the Shapley effects with green horizontal segments.}
    \label{boxplot}
\end{figure}

Here, we carry out a numerical study based on the following function:
$$
f(x)=\cos(z)+z-100+0.2\sin(10z),\;\;z=\sum_{i=1}^{p} \left(1+\frac{i-1}{p-1}\right) x_i^2,
$$
with a number of inputs $p=15$. For simplicity, we take the Gaussian input distribution $\mathcal{N}(0,\Gamma)$ where the covariance matrix $\Gamma$ is block-diagonal, with 3 blocks of size 5.

In Figure \ref{boxplot}, we present the estimates of the Shapley effects achieved with the two approaches for the same computation time (140 seconds for computing the meta-model and estimating the Shapley effects of this meta-model). Each method is run 1000 times to obtain the boxplots. At each iteration, we generate a new sample $(X_i,f(X_i))_{i\leq N_{samp}}$ and we compute a new model $\tilde{f}_{tot}$ or $\tilde{f}_{add}$. We also represent the true values of the Shapley effects, computed by taking the mean of 1000 estimates given by shapleyPermRand with the true model $f$ (and the parameters $Nv=10^6$, $m = 10^4$).

We remark that the estimates obtained from the additive model $\tilde{f}_{add}$ and from Corollary \ref{coroshapley} are more accurate than the estimates with the total model $\tilde{f}_{tot}$ for the same computation time.

%In conclusion, to estimate the Shapley effects of a costly numerical code with independent groups of inputs, we suggest to create a block-additive meta-model and to estimate the Shapley effects of this meta-model using Corollary \ref{coroshapley}.
In conclusion, the use of a block-additive model, and the estimation of the Shapley effects for this meta-model using corollary \ref{coroshapley} appears to be relevant to estimate the Shapley effects of a costly numerical code with independent groups of input variables.
As we have shown, this can results in an improved accuracy, even for codes that are not block-additive.

\section{Conclusion}
In this paper, we give new theoretical results about the variance-based sensitivity indices for independent groups of inputs. These results drastically reduce the computational cost of these indices when the model is block-additive. Then, we apply these results to the linear Gaussian framework and we suggest two algorithms: the first one for the general case and the second one for a block diagonal covariance matrix. These algorithms compute efficiently the theoretical values of the variance-based sensitivity indices. Numerical experiments on Shapley effects computations highlight this efficiency and the benefit compared to existing methods. We also suggest a method based on our results to improve the estimates of the Shapley effects when we have independent groups of variables, and a general model that is not block-additive.

\section*{Acknowledgements}
We are very grateful to Bertrand Iooss for his ideas and his advices.
We also would like to thank Pietro Mosca and Laura Clouvel for their contribution to the implementation of the industrial application, and Arnaud Grivet Sebert for the constructive discussions that helped to improve the paper. The authors are grateful to the referees for their constructive suggestions.

\section*{Appendix}

\textbf{Proof of Proposition \ref{propcentrale}:}\\
We use the Hoeffding decomposition of $g$:
\begin{eqnarray*}
V_u&=&\V(\E(Y|X_u))\\
&=&\V\left[\sum_{w\subseteq[1:k]}\E(g_w(A_w)|X_{u\cap C_w})\right]\\
&=&\sum_{w\subseteq[1:k]}\V\left[\E(g_w(A_w)|X_{u\cap C_w})\right]\\
&=&\sum_{w\subseteq[1:k]} V_{u\cap C_w}^{g,w}.
\end{eqnarray*}\qed

\bigskip

\textbf{Proof of Proposition \ref{propSobol}:}\\
We have
\begin{eqnarray*}
S_u&=&\frac{1}{\V(Y)}\sum_{v\subseteq u} (-1)^{|u|-|v|}V_v\\
&=&\frac{1}{\V(Y)}\sum_{v\subseteq u} (-1)^{|u|-|v|}\sum_{w\subseteq[1:k]} V_{v\cap C_w}^{g,w}\\
&=&\frac{1}{\V(Y)}\sum_{w\subseteq[1:k]} \sum_{v\subseteq u} (-1)^{|u|-|v|}V_{v\cap C_w}^{g,w}\\
&=&\frac{1}{\V(Y)}\sum_{w\subseteq[1:k]} \sum_{v_1\subseteq u\cap C_w} V_{v_1}^{g,w} \displaystyle\sum_{\substack{v_2 \subseteq -C_w,\\ \text{s.t.} v_1 \cup v_2 \subseteq u}}(-1)^{|u|-|v_1|-|v_2|}
\end{eqnarray*}
One can remark that, if $u\setminus C_w\neq \emptyset$,
$$
\sum_{\substack{v_2 \subseteq -C_w,\\ \text{s.t.} v_1 \cup v_2 \subseteq u}}(-1)^{|u|-|v_1|-|v_2|} =(-1)^{|u|-|v_1|}\sum_{n=0}^{|u\setminus C_w|}\begin{pmatrix}
|u\setminus C_w|
\end{pmatrix}
(-1)^{n}=0.
$$
Thus,
\begin{eqnarray*}
S_u&=&\frac{1}{\V(Y)}\displaystyle\sum_{\substack{w\subseteq[1:k],\\ u\subseteq C_w}} \sum_{v_1\subseteq u} (-1)^{|u|-|v|} V_{v_1}^{g,w} \\
&=&\displaystyle\sum_{\substack{w\subseteq[1:k],\\ u\subseteq C_w}} \frac{\V(g_w(A_w))}{\V(Y)}S_u^{g,w}.
\end{eqnarray*}
\qed

\bigskip

\textbf{Proof of Proposition \ref{propshapley}:}\\

\begin{eqnarray*}
\eta_i&=& \frac{1}{p\V(Y)} \sum_{u\subseteq -i} \begin{pmatrix} p-1\\ |u| \end{pmatrix}^{-1} (V_{u\cup \{i\} }-V_u)\\
&=& \frac{1}{p\V(Y)} \sum_{u\subseteq -i} \begin{pmatrix} p-1\\ |u| \end{pmatrix}^{-1} \sum_{w\subseteq[1:k]}(V_{(u \cup i) \cap C_w}^{g,w}-V_{u\cap C_w}^{g,w})\\
&=& \frac{1}{p\V(Y)} \displaystyle\sum_{\substack{w\subseteq[1:k],\\ \text{s.t. }j(i)\in w}} \sum_{u\subseteq -i} \begin{pmatrix} p-1\\ |u| \end{pmatrix}^{-1} (V_{(u \cup i) \cap C_w}^{g,w}-V_{u\cap C_w}^{g,w})\\
&=& \frac{1}{p\V(Y)} \displaystyle\sum_{\substack{w\subseteq[1:k],\\ \text{s.t. } j(i)\in w}} \sum_{u\subseteq C_w\setminus \{i\}}  \left[\sum_{v \subseteq -C_w}  \begin{pmatrix} p-1\\ |u\cup v| \end{pmatrix}^{-1}  \right](V_{u\cup i}^{g,w}-V_{u}^{g,w})\\
&=& \frac{1}{p\V(Y)}\displaystyle\sum_{\substack{w\subseteq[1:k],\\ \text{s.t. } j(i)\in w}} \sum_{u\subseteq C_w\setminus \{i\}} \left[\sum_{j=0}^{p-|C_{w}|}\begin{pmatrix}
p-|C_{w}|\\ j
\end{pmatrix}\begin{pmatrix}
p-1 \\ |u|+j
\end{pmatrix}^{-1}\right](V_{u\cup i}^{g,w}-V_{u}^{g,w}).
\end{eqnarray*}
It remains to prove the following equation:
\begin{eqnarray}\label{binom}
\frac{1}{p}\sum_{j=0}^{p-|C_{w}|}\begin{pmatrix}
p-|C_{w}|\\ j
\end{pmatrix}\begin{pmatrix}
p-1 \\ |u|+j
\end{pmatrix}^{-1}=\frac{1}{|C_{w}|}\begin{pmatrix}
|C_{w}|-1 \\ |u|
\end{pmatrix}^{-1}.
\end{eqnarray}
In the interest of simplifying notation, until the end of the proof, we will write $u$ (resp. $c$) instead of $|u|$ (resp. $|C_{w}|$).
We can verify that the equation \eqref{binom} is equivalent to the following equations:
\begin{eqnarray}\label{bino}
&&\frac{1}{p}\sum_{j=0}^{p-c}\frac{(p-c)!}{j!(p-c-j)!} \frac{(u+j)!(p-1-u-j)!}{(p-1)!}=\frac{1}{c}\frac{u!(c-1-u)!}{(c-1)!} \nonumber \\
& \Longleftrightarrow & \sum_{j=0}^{p-c}\begin{pmatrix}
u+j\\u
\end{pmatrix} \begin{pmatrix}
p-1-u-j\\c-u-1
\end{pmatrix}=\begin{pmatrix}
p\\ c
\end{pmatrix}.
\end{eqnarray}
We will show \eqref{bino}. Now, we can remark that we have:
\begin{equation}
\frac{x^c}{(1-x)^{c+1}}=x\frac{x^u}{(1-x)^{u+1}}\frac{x^{c-u-1}}{(1-x)^{c-u}}.
\end{equation}
Giving their power series, we have:
\begin{equation}\label{power}
x \left( \sum_{k\geq0}\begin{pmatrix}k\\ u\end{pmatrix}x^k \right) \left( \sum_{k'\geq0}\begin{pmatrix}k'\\ c- u-1 \end{pmatrix}x^{k'} \right)=  \sum_{k''\geq0}\begin{pmatrix} k''\\ c\end{pmatrix}x^{k''}.
\end{equation}
We have the equality of the coefficient of $x^p$. Then:
\begin{eqnarray*}
\begin{pmatrix}
p\\c
\end{pmatrix} &=&
\sum_{k+k'=p-1}\begin{pmatrix}
k\\u
\end{pmatrix}\begin{pmatrix}
k'\\c-u-1
\end{pmatrix}
=\sum_{k=u}^{p-1}\begin{pmatrix}
k\\u
\end{pmatrix}\begin{pmatrix}
p-1-k\\c-u-1
\end{pmatrix}
=\sum_{j=0}^{p-1-u}\begin{pmatrix}
u+j\\u
\end{pmatrix}\begin{pmatrix}
p-1-u-j\\c-u-1
\end{pmatrix}\\
&=&\sum_{j=0}^{p-c}\begin{pmatrix}
u+j\\u
\end{pmatrix}\begin{pmatrix}
p-1-u-j\\c-u-1
\end{pmatrix}.
\end{eqnarray*}
For the last equality, we remark that if $j>p-c$, then $p-1-u-j<c-u-1$ so the last terms of the sum are equal to zero. We have proven \eqref{bino}. To conclude, we have
\begin{eqnarray*}
\eta_i&=&\frac{1}{p\V(Y)}\displaystyle\sum_{\substack{w\subseteq[1:k],\\ \text{s.t. } j(i)\in w}} \sum_{u\subseteq C_w\setminus \{i\}} \left[\sum_{j=0}^{p-|C_{w}|}\begin{pmatrix}
p-|C_{w}|\\ j
\end{pmatrix}\begin{pmatrix}
p-1 \\ |u|+j
\end{pmatrix}^{-1}\right] (V_{u \cup \{i\}} -V_{u})\\
&=&\displaystyle\sum_{\substack{w\subseteq[1:k],\\ \text{s.t. } j(i)\in w}}S_w^g  \frac{1}{|C_{w}| \V(g_w(A_w)}\sum_{u\subseteq C_w\setminus \{i\}} \begin{pmatrix}
|C_{w}|-1 \\ |u|
\end{pmatrix}^{-1}(V_{u\cup i}^{g,w}-V_{u}^{g,w})\\
&=& \displaystyle\sum_{\substack{w\subseteq[1:k],\\ \text{s.t. } j(i)\in w}}S_w^g  \eta_i^{g,w}.
\end{eqnarray*}
\qed

\section*{References}
\bibliographystyle{alpha}

\bibliography{bibfile}

\end{document}